\documentclass[12pt]{amsart}
\usepackage{amscd,amssymb,amsmath}
\usepackage{graphicx}

\newtheorem{thm}[subsection]{Theorem}
\newtheorem{lemma}[subsection]{Lemma}
\newtheorem{prop}[subsection]{Proposition}
\newtheorem{cor}[subsection]{Corollary}

\newtheorem{rk}[subsection]{Remark}
\newtheorem{defn}[subsection]{Definition}

\numberwithin{equation}{section} 

\newcommand{\bea}{\begin{eqnarray}}
\newcommand{\eea}{\end{eqnarray}}

\newcommand{\Z}{\mathbb{Z}}
\newcommand{\Q}{\mathbb{Q}}

\newcommand{\N}{\mathbb{N}}

\begin{document}
\title [Equation $x^q=a$ in field of $p$-adic numbers]
{Solvability criteria for the equation $x^q=a$ in the field of $p$-adic numbers}

\author {J.M. Casas, B.A.Omirov, U.A. Rozikov}

 \address{J.\ M.\ Casas\\ Department of Applied Mathematics, E.U.I.T. Forestal,  University of Vigo, 36005, Pontevedra, Spain.}
 \email {jmcasas@uvigo.es}
  \address{B.\ A.\ Omirov\\ Institute of mathematics and information technologies,
Tashkent, Uzbekistan.} \email {omirovb@mail.ru}
 \address{U.\ A.\ Rozikov\\ Institute of mathematics and information technologies,
Tashkent, Uzbekistan.} \email {rozikovu@yandex.ru}

\begin{abstract} We establish the solvability criteria for the equation $x^q=a$ in the field of $p$-adic numbers, for any $q$ in two cases:
 (i) $q$ is not divisible by $p$; (ii) $q=p$. Using these criteria we show that any $p$-adic number can be represented
 in finitely many different forms and we describe the algorithms to obtain the corresponding representations. Moreover it is showed that solvability problem of $x^q=a$ for any $q$ can be reduced to
 the cases (i) and (ii).

{\it AMS classifications (2010):} 11Sxx

{\it Keywords:} $p$-adic number; solvability of an equation; congruence.
\end{abstract}
\maketitle
\section{Introduction} \label{sec:intro}

The $p$-adic number system for any prime number $p$ extends the ordinary arithmetic of the rational numbers in
a way different from the extension of the rational number system to the real and complex number systems.
This extension is achieved by an alternative interpretation of the concept of absolute value.

First described by Kurt Hensel in 1897, the $p$-adic numbers were motivated primarily by an attempt to
bring the ideas and techniques of power series methods into number theory.
Their influence now extends far beyond this. For example, the field of $p$-adic analysis essentially
provides an alternative form of calculus.

More formally, for a given prime $p$, the field $\Q_p$ of $p$-adic numbers is a completion of the rational numbers.
On the field $\Q_p$ is also given a topology derived from a metric, which is itself derived from an alternative
valuation on the rational numbers. This metric space is complete in the sense that every Cauchy sequence converges
to a point in $\Q_p$. This is what allows the development of calculus on $\Q_p$ and it is the interaction
of this analytic and algebraic structure which gives the $p$-adic number systems their power and utility.

 For about a century after the discovery of $p$-adic numbers, they were
mainly considered as objects of pure mathematics. However,
numerous applications of these numbers to theoretical physics have
been proposed in papers \cite{5,2,3,4}, to quantum mechanics
\cite{6}, to $p$-adic - valued physical observables \cite{6} and
many others \cite{K,V}.  As in real case, to solve a problem in
$p$-adic case there arises an equation which must be solved in the
field of $p$-adic numbers (see for example \cite{M1, M2, M3, V}).
For classification problems of varieties of algebras over a field
$\Q_p$ of $p$-adic numbers one has to solve an equation of the
form: $x^2=a$. It is well known the criteria of solvability of
this equation. In fact, in the classification of Leibniz algebras
of dimensions less than 4 over a field $\Q_p$ (see \cite{Ay},
\cite{Abror}) it is enough to solve the equation $x^2=a$. However,
in the classifications tasks of larger dimensions one has to solve
an equation of the form $x^q=a$, $q\geq 2$, in $\Q_p$. In this
paper we present a criteria for solvability of the equation
$x^q=a$ in $\Q_p$ (where $p$ is a fixed prime number) for
arbitrary $q$ in two cases: $(q,p)=1$ and $q=p$. Moreover, in the cases of existing the solutions of the equation we present the algorithm of finding the solutions.
Also we show that any equation $x^q=a$ in $Q_p$ can be reduced to both cases. Note that in \cite{MS} the same criterion has been proved
by a different method.

\section{Preliminaries}

\subsection{Solvability of congruences.} The sources of the information in this subsection are \cite{A,B,N}.  If $n$ is a positive integer, the integers between 1 and $n-1$ which are coprime to $n$ (or equivalently, the congruence classes coprime to $n$) form a group with multiplication modulo $n$ as the operation; it is denoted by $\Z^\times_n$ and is called the {\it group of units modulo} $n$ or the group of primitive classes modulo $n$. Multiplicative group of integers modulo $n$ is cyclic if and only if $n$ is equal to $1, 2, 4, p^k$ or $2 p^k$, where $p^k$ is a power of an odd prime number. A generator of this cyclic group is called {\it a primitive root modulo} $n$ or a primitive element of $\Z^\times_n$.

 For any integers $a, b$ and $n$, we say that $a$ is congruent to $b$ modulo $n$ (notated $a \equiv b \mod n$)
if $n\mid (b-a)$. The order of $\Z^\times_n$ is given by Euler's totient function  $\varphi(n)$. Euler's theorem says that $a^{\varphi(n)}\equiv 1 \mod n$ for every $a$ coprime to $n$; the lowest power of $a$ which is $\equiv 1 \mod n$ is called the multiplicative order of $a$ modulo $n$. In other words, for $a$ to be a primitive root modulo $n$, $\varphi(n)$ has to be the smallest power of $a$ which is
 congruent to 1  $\mod n$.
\begin{lemma}\label{l1} Suppose that $m\in \N$ has a primitive root $r$.
If $a$ is a positive integer with $(a,m) = 1$, then there is a unique
integer $x$ with $1\leq x\leq \varphi(m)$ such that
$$r^x\equiv a \mod m.$$
\end{lemma}

\begin{defn} If $m\in \N$ has a primitive root $r$ and $\alpha$ is a
positive integer with $(\alpha,m) = 1$, then the unique integer $x$,
$1 \leq x \leq \varphi(m)$ and $r^x\equiv \alpha \mod m$ is called the index (or discrete
logarithm) of $\alpha$ to the base $r$ modulo $m$ and denoted by ${\rm ind}_r\alpha$.
\end{defn}
In particular, $r^{{\rm ind}_ra}\equiv a \mod m.$

\begin{thm}\label{t1} Let $m$ be a positive integer with primitive root
$r$. If $a, b$ are positive integers coprime to $m$ and $k$ is a positive integer, then

(i) ${\rm ind}_r1\equiv 0 \mod \varphi(m)$

(ii) ${\rm ind}_r(ab)\equiv {\rm ind}_ra + {\rm ind}_rb \mod\varphi(m)$

(iii) ${\rm ind}_ra^k\equiv k\cdot {\rm ind}_ra \mod \varphi(m)$.
\end{thm}

\begin{thm}\label{tv} If $p$ is a prime number, $\alpha\in \N$, $m$ is equal to $p^\alpha$ or $2p^{\alpha}$, $(n,\varphi(m))=d$ then the congruence
$$x^n\equiv a\mod m$$  has solution if and only if $d$ divides ${\rm ind}_x a$. In case of solvability the congruence equation has $d$ solutions.
\end{thm}
{\it Fermat's Little Theorem} says: Let $a$ be a nonzero integer, and let $p \nmid a$
be prime. Then $a^{p-1} \equiv 1 \mod p$.

\begin{prop}\label{p1} Let $a, b, n \in \Z$ with $n\ne 0$. The congruence
$ax \equiv b \mod n$
has solutions if and only if $(a, n)\mid b$. When the congruence has a solution $x_0\in  \Z$
then the full solution set is $ \{{x_0 + tn\over (a, n)} : t \in \Z\}$.
\end{prop}
It follows that the equation $ax = b$ in $\Z_n$ has $(a, n)$ solutions. In particular, if
$(a, n) = 1$, then the equation $ax = b$ has  unique solution in $\Z_n$.

\subsection{Divisibility of binomial coefficients.} (see \cite{H}) In 1852, Kummer proved that if $m$ and $n$ are nonnegative integers and $p$ is a prime number,
then the largest power of $p$ dividing ${m+n\choose m}$ equals $p^c$, where $c$ is the number of carries when $m$ and $n$ are added in base $p$. Equivalently, the exponent of a prime $p$ in ${n\choose k}$  equals the number of nonnegative integers $j$ such that the fractional part of $k/p^j$ is greater than the fractional part of $n/p^j$. It can be deduced from this that ${n\choose k}$  is divisible by $n/gcd(n,k)$. Another fact: an integer $n\geq 2$ is prime if and only if all the intermediate binomial coefficients
${n\choose k}$, $k=1,\dots,n-1$ are divisible by $n$.

\subsection{ $p$-adic numbers.} Let $\Q$ be the field of rational numbers. Every rational number $x\ne 0$ can be represented
in the form $x = p^r{n\over m}$, where $r, n\in \Z$, $m$ is a positive integer, $(p, n) = 1$, $(p, m) = 1$ and $p$
 is a fixed prime number. The p-adic norm of $x$ is given by
$$|x|_p=\left\{\begin{array}{ll}
p^{-r}\ \ \mbox{for} \ \ x\ne 0\\
0\ \ \mbox{for} \ \ x = 0.
\end{array}\right.
$$
This norm satisfies the so called strong triangle inequality
$$|x+y|_p\leq \max\{|x|_p,|y|_p\},$$
and this is a non-Archimedean norm.

The completion of $\Q$ with respect to $p$-adic norm defines the $p$-adic field
which is denoted by $\Q_p$. Any $p$-adic number $x\ne 0$ can be uniquely represented
in the canonical form
\begin{equation}\label{ek}
x = p^{\gamma(x)}(x_0+x_1p+x_2p^2+\dots),
\end{equation}
where $\gamma=\gamma(x)\in \Z$ and $x_j$ are integers, $0\leq x_j \leq p - 1$, $x_0 > 0$, $j = 0, 1,2, ...$ (see
more detail \cite{K,V,S}). In this case $|x|_p = p^{-\gamma(x)}$.

\begin{thm}\label{tx2} \cite{B}, \cite{K}, \cite{V}  In order to the equation
$x^2 = a$, $0\ne a =p^{\gamma(a)}(a_0 + a_1p + ...), 0\leq a_j \leq p - 1$, $a_0 > 0$
has a solution $x\in \Q_p$, it is necessary and sufficient that the following conditions
are fulfilled:

i) $\gamma(a)$ is even;

ii) $a_0$ is a quadratic residue modulo $p$ if $p\ne 2$; $a_1 = a_2 = 0$ if $p = 2$.
\end{thm}

In this paper we shall generalize this theorem.

\section{Equation $x^q=a$.}

In this section we consider the equation $x^q=a$ in $\Q_p$, where $p$ is a fixed prime number, $q\in \N$ and $a\in \Q_p$. Our goal is to find conditions under which the equation has a solution $x\in \Q_p$. The case $q=2$ is well known (see Theorem \ref{tx2}), therefore we consider $q>2$.

We need the following

\begin{lemma}\label{l2} The following are true

(i)
\begin{equation}\label{e3}
\left(\sum_{i=0}^\infty x_ip^i\right)^q=x_0^q+\sum_{k=1}^\infty\left(qx_0^{q-1}x_k+N_k(x_0,x_1,\dots,x_{k-1})\right)p^k,\end{equation}
where $x_0\ne 0$, $0\leq x_j\leq p-1$, $N_1=0$ and for $k\geq 2$
\begin{equation}\label{e4}
N_k=N_k(x_0,\dots,x_{k-1})=\\ \sum_{{m_0,m_1,\dots,m_{k-1}:\atop
\sum_{i=0}^{k-1}m_i=q,\ \
\sum_{i=1}^{k-1}im_i=k}} {q!\over m_0!m_1!\dots m_{k-1}!}x_0^{m_0}x_1^{m_1}\dots x_{k-1}^{m_{k-1}}.
\end{equation}

(ii) Let $q=p$ be a prime number,
then $p\mid N_k$ if and only if $p\nmid k$.

\end{lemma}
\proof (i) Formulas (\ref{e3}) and (\ref{e4}) easily can be obtained by using multinomial theorem.

(ii) It is known the following formula:
\begin{equation}\label{e5}
{p!\over m_0!m_1!\dots m_{k-1}!}={m_0 \choose m_0}{m_0+m_1 \choose m_1}\dots {m_0+m_1+\dots+m_{k-1} \choose m_{k-1}}.
\end{equation}

By condition $\sum_{i=0}^{k-1}m_i=p$ we get $0\leq m_i\leq p$. Moreover if $m_{i_0}=p$ for some $i_0$, then $m_{i}=0$ for all $i\ne i_0$. In this case from the  condition $\sum_{i=1}^{k-1}im_i=k$ we obtain $k=i_0p$.

Assume  $p\mid k$, i.e., $k=pk_0$ for some $k_0$, where  $1\leq k_0<k-1$. Putting $m_{k_0}=p$ and $m_i=0$ for $i\ne k_0$, we get $N_k(x_0,\dots,x_{k-1})=x_{k_0}^p+S_k,$ where
$$S_k=\sum_{{0\leq m_0,m_1,\dots,m_{k-1}<p:\atop
\sum_{i=0}^{k-1}m_i=p,\ \ \sum_{i=1}^{k-1}im_i=k}} {p!\over m_0!m_1!\dots m_{k-1}!}x_0^{m_0}x_1^{m_1}\dots x_{k-1}^{m_{k-1}}.$$
Since $x_{k_0}^p$ is not a multiple of $p$ (otherwise, if we assume that $r^p=pT$ for some $r\in \{2,\dots,p-1\}$ and $T\in \N$, then
it follows that $r$ divides $T$ since $p$ is a prime number. Hence, $r^{p-1}=pT',$ where $T'=T/r$. By the same way we get $r^{p-2}=pT''$ and iterating the process finally we get
 $1=p\tilde{T}$, for some $\tilde {T}\in \N$, which is not possible), but $S_k$ is divisible by $p$
 since each coefficient of
$S_k$ contains (see formula (\ref{e5})) a factor
 ${p\choose m_{i_0}}$ (with $0<m_{i_0}<p$) which is divisible by $p$ thanks to the divisibility property of binomial coefficients  mentioned in the previous section. Therefore
 $p\nmid N_k$.

Assume  $p\nmid k$ then by the arguments mentioned above we get $m_i<p$, for all $i=0,\dots,k-1$,
therefore each term of $N_k$ is divisible by $p$ consequently    $p\mid N_k$.
\endproof

3.1. {\sl THE CASE $(q,p)=1$.} In this subsection we are going to
analyze under what conditions the equation $x^q=a$ has solution in
$Q_p$, when $q$ and $p$ are coprimes.  In this case the following
is true.

\begin{thm}\label{t2} Let $q>2$ and $(q,p)=1$. The equation
\begin{equation}\label{e1}
x^q=a,
\end{equation}
 $0\ne a=p^{\gamma(a)}(a_0+a_1p+\dots)$, $0\leq a_j\leq p-1$, $a_0\ne 0$, has a solution $x\in \Q_p$ if and only if

1) $q$ divides $\gamma(a)$;

2) $a_0$ is a $q$ residue$\mod p$.
\end{thm}

\proof {\it Necessity}.  Assume that equation (\ref{e1}) has a solution
$$ x=p^{\gamma(x)}\left(x_0+x_1p+\dots\right), \ \ 0\leq x_j\leq p-1, x_0\ne 0,$$ then
\begin{equation}\label{e2}
p^{q\gamma(x)}\left(x_0+x_1p+\dots\right)^q=p^{\gamma(a)}(a_0+a_1p+\dots).
\end{equation}
Consequently, $\gamma(a)=q\gamma(x)$ and $a_0\equiv x_0^q \mod p$.

{\it Sufficiency.} Let $a$ satisfies the conditions 1) and 2). We construct a solution $x$ of equation (\ref{e1}) using the idea of reduction to canonical form of a $p$-adic number using a system of carries.
Put
$$\gamma(x)={1\over q}\gamma(a).$$

Then by condition 2) and due to $1\leq a_0\leq p-1$ there exists $x_0$ such that
\begin{equation*}\label{e7}
x_0^q\equiv a_0 \mod p, \ \ 1\leq x_0\leq p-1.
\end{equation*}
In other words, there exists $M_1(x_0)$ such that $x_0^q=a_0+M_1(x_0)p.$

Using the notations of Lemma \ref{l2} due to the fact that $qx_0^{q-1}$ is not a multiple of $p,$ there exists $x_1$ such that
\begin{equation*}\label{e71}
qx_0^{q-1}x_1+N_1(x_0)+M_1(x_0)\equiv a_1 \mod p, \ \ 1\leq x_1\leq p-1.
\end{equation*}
Therefore, there exists $M_2(x_0,x_1)$ such that $qx_0^{q-1}x_1+N_1(x_0)+M_1(x_0)=a_1+M_2(x_0,x_1)p.$

Proceeding in this way, we find the existence of $x_n$ such that
$$
qx_0^{q-1}x_n+N_n(x_0,\dots,x_{n-1})+M_n(x_0,\dots,x_{n-1})\equiv a_n \mod p, \ \ 1\leq x_n\leq p-1.
$$
and $M_{n+1}(x_0,\dots,x_n)$ such that
\begin{equation}\label{e72}
qx_0^{q-1}x_n+N_n(x_0,\dots,x_{n-1})+M_n(x_0,\dots,x_{n-1})=a_n+M_{n+1}(x_0,\dots, x_{n})p
\end{equation}
for any $n\in \mathbb{N}.$

Now from Lemma \ref{l2} and equality (\ref{e72}) it follows that
$$\left(\sum_{i=0}^\infty x_ip^i\right)^q=x_0^q+\sum_{k=1}^\infty\left(qx_0^{q-1}x_k+N_k(x_0,x_1,\dots,x_{k-1})\right)p^k=$$
$$=a_0+M_1(x_0)p+\sum_{k=1}^\infty\left(a_k-M_k(x_0,\dots,x_{k-1})+M_{k+1}(x_0,\dots,x_k)p\right)p^k=a_0+\sum_{k=1}^\infty a_k p^k.$$
Hence, we found solution $x=\displaystyle\sum_{i=0}^\infty x_ip^i$ of equation (\ref{e1}) in its canonical form.\endproof

\begin{rk} The condition 2) of Theorem \ref{t2} is always satisfied if $p=2$, consequently $x^q=a$ has a solution in $\Q_2$ for any odd $q$
and any $a\in \Q_2$ with $\gamma(a)$ divisible by $q$.
\end{rk}

\begin{cor}\label{c1} Let $q$ be a prime number such that $q<p$ and $\eta$ be a unity (i.e. $|\eta|_p=1$) which is not $q$-th power of some $p$-adic number. Then $p^i\eta^j$, $i,j=0,1,\dots,q-1$ ($i+j\ne 0$) is not $q$-th power of some $p$-adic number.
\end{cor}
\proof We shall check that the conditions of Theorem \ref{t2} fail under the hypothesis.
It is easy to see that $\gamma(p^i\eta^j)=i$ for all $i=1,\dots,q-1$ and $j=0,\dots,q-1$, consequently $q$ does not divide $\gamma(p^i\eta^j)$, i.e. condition 1) is not satisfied. But for $\eta^j$, $i=0$, $j=2,\dots, q-1$ the condition $1)$ is satisfied, therefore we shall check condition $2)$.

Consider the decomposition $\eta=\eta_0+\eta_1p+\dots$, then $\eta^j=\eta_0^j+j\eta_0^{j-1}\eta_1p+\dots$. It is known (see Theorem \ref{t1}) that $\eta^j_0\equiv a_0^q \mod p$ has a solution $\eta_0$ if and only if
${\rm ind}_{a_0}\eta^j_0$ is divisible by $d=(p-1,q)$.

Since $\eta$ is a unity  which is not $q$-th power of some $p$-adic number we have $\eta_0\equiv a_0^q \mod p$ has not solution thus ${\rm ind}_{a_0}\eta_0$ is not divisible by $d=(p-1,q)$. This property
implies that for a prime $q$ with $q<p$ the prime number $p$ has a form $p=qk+1$. Consequently, $d=(p-1,q)=(qk,q)=q$.

If ${\rm ind}_{a_0}\eta_0=dl+r$, then ${\rm ind}_{a_0}\eta_0^j\equiv j(ql+r)\mod(p-1)$, i.e., ${\rm ind}_{a_0}\eta_0^j= j(ql+r)+M(p-1)=q(jl+Mk)+jr$, but
 since $0<r<q-1$, $2\leq j\leq q-1$ and $q$ is a prime number, $jr$ is not divisible by $q$. Thus ${\rm ind}_{a_0}\eta_0^j$ is not divisible by $d=q$, hence condition $2)$ is not satisfied.
 \endproof

\begin{cor}\label{c2} Let $q$ be a prime number such that $q<p=qk+1$, for some $k\in \N$ and $\eta$ be a unity which is not $q$-th power of some $p$-adic number. Then any $p$-adic number $x$ can be written in one of the following forms $x=\varepsilon_{ij}y_{ij}^q$, where $\varepsilon_{ij}\in \{p^i\eta^j: i,j=0,1,\dots,q-1\}$ and $y_{ij}\in \Q_p$.
\end{cor}
\proof Let $\eta=\eta_0+\eta_1p+\dots$ and $\mu=\mu_0+\mu_1p+\dots$ be unities which are not $q$-th power of some $p$-adic numbers.

We shall show that there exists $i\in \{1,\dots,q-1\}$ such that $\mu=\eta^iy^q$ for some $y\in \Q_p$.

Note that $\eta^i$ and ${1\over \eta^i}$, $i=1,\dots,q-1$ are not $q$-th power of some $p$-adic numbers.
Consider $\mu=\mu_0+\mu_1p+\dots$ and ${1\over\eta^i}=c_{0i}+c_{1i}p+\dots$, then
${\mu\over \eta^i}=\mu_0c_{0i}+(\mu_1c_{0i}+\mu_0c_{1i})p+\dots$

It is easy to see that $\gamma(\mu/\eta^i)=0$, consequently condition $1)$ of Theorem \ref{t2} is satisfied. Indeed, if $\mu_0c_{0i}=pM$, then since $p$ is a prime number, the last equality is possible if and only if $\mu_0$ and $c_{0i}$ divide $M$, consequently we get $1=p\tilde{M}$, which is impossible.

Now we shall check condition $2)$ of Theorem \ref{t2}.

The equation $x_0^q\equiv \mu_0c_{0i}\mod p$ has a solution if and only if ${\rm ind}_{x_0}\mu_0c_{0i}$ is divisible by $q=(q,p-1)$. We have that $c_{0i}\equiv c_{01}^i\mod p$, consequently  ${\rm ind}_{x_0}\mu_0c_{0i}\equiv {\rm ind}_{x_0}\mu_0+i \,{\rm ind}_{x_0}c_{01} \mod (p-1)$. Assume ${\rm ind}_{x_0}\mu_0=s$, ${\rm ind}_{x_0}c_{01}=r$, $s,r=1,\dots, q-1$ and take $i$ such that
$ir\equiv q-s \mod (p-1)$ (the existence of $i$ follows from the fact that $(r,p-1)=1$). It is clear that if $s$ varies from 1 to $q-1$ then $i$ also varies in $\{1,\dots,q-1\}$. Consequently we get ${\rm ind}_{x_0}\mu_0c_{0i}\equiv s+ir\mod (p-1)\equiv q \mod (p-1)$. Hence the condition $2)$
in Theorem \ref{t2} is also satisfied, so there is an $i$ such that $\mu=\eta^i y_i^q$, for some $y_i\in \Q_p$.

For $x\in \Q_p$, let $x=p^{\gamma(x)}(x_0+x_1p+\dots)$ denote $\nu=x_0+x_1p+\dots$ If $\nu$ satisfies conditions of Theorem \ref{t2}, then
$\nu=y^q$ and $x=p^{\gamma(x)}y^q$. If $\nu$ does not satisfy conditions of Theorem \ref{t2}, then (as it was showed above) there exists $i$ such that
$\nu=\eta^iy^q$ and in this case we get $x=p^{\gamma(x)}\eta^iy^q$. Taking $\gamma(x)=qN+j$ completes the proof.
\endproof
\begin{cor}\label{c3a} Let $q$ be a prime number such that $q<p\ne qk+1$, for any $k\in \N$. Then any $p$-adic number $x$ can be written in one of the following forms $x=\varepsilon_{i}y_{i}^q$, where $\varepsilon_{i}\in \{p^i: i=0,1,\dots,q-1\}$ and $y_{i}\in \Q_p$.
\end{cor}
\proof Using arguments in the proof of Corollary \ref{c2} we conclude that if $p\ne qk+1$ then any unity is $q$-th power of a $p$-adic number. Hence for $x\in \Q_p$ with $\gamma(x)=qN+i$ we get $x=p^iy^q$, $y\in \Q_p$.
\endproof

{3.2. \sl THE CASE $q=p$.}  Now we are going to analyze the solvability conditions for the equation $x^p=a$  in $Q_p$. In this case, the following is true.

\begin{thm}\label{t3} Let $q=p$. The equation $x^q=a$,
 $0\ne a=p^{\gamma(a)}(a_0+a_1p+\dots)$, $0\leq a_j\leq p-1$, $a_0\ne 0$, has a solution $x\in \Q_p$ if and only if

 (i) $p$ divides $\gamma(a)$;

 (ii) $a_0^p\equiv a_0+a_1p\mod p^2$.
\end{thm}

\proof {\it Necessity}.  Assume that equation (\ref{e1}) has a solution
$$ x=p^{\gamma(x)}\left(x_0+x_1p+\dots\right), \ \ 0\leq x_j\leq p-1, x_0\ne 0,$$ then using Lemma \ref{l2} we get
\begin{equation}\label{e8}
a=p^{p\gamma(x)}\left(x_0+x_1p+\dots\right)^p=p^{p\gamma(x)}\left(x_0^p+\sum_{k=1}^\infty(px_0^{p-1}x_k+N_k)p^k\right).
\end{equation}
Using part $(ii)$ of Lemma \ref{l2}, we get
$${\rm RHS \ \ of \ \ (\ref{e8})}=p^{p\gamma(x)}\left(x_0^p+\sum_{{k=1\atop p\mid k}}^\infty x_0^{p-1}x_kp^{k+1}+
\sum_{{k=1\atop p\mid k}}^\infty N_kp^k+\right.$$ $$\left.\sum_{{k=1\atop p\nmid k}}^\infty(x_0^{p-1}x_k+p^{-1}N_k)p^{k+1}\right)=$$
$$p^{p\gamma(x)}\left(x_0^p+\sum_{k=1}^\infty x_0^{p-1}x_{pk}p^{pk+1}+
\sum_{k=1}^\infty N_{pk}p^{pk}+\right.$$ $$\left.\sum_{i=1}^{p-1}\sum_{k=0}^\infty(x_0^{p-1}x_{pk+i}+p^{-1}N_{pk+i})p^{pk+i+1}\right)=$$
$$p^{p\gamma(x)}\left(x_0^p+\sum_{i=1}^{p-2}(x_0^{p-1}x_i+p^{-1}N_i)p^{i+1}+\right.$$ $$\left.\sum_{k=1}^\infty (x_0^{p-1}x_{pk-1}+N_{pk}+p^{-1}N_{pk-1})p^{pk}
+\right.$$
\begin{equation}\label{e9}
\left.\sum_{k=1}^\infty x_0^{p-1}x_{pk}p^{pk+1}+\sum_{i=1}^{p-2}\sum_{k=1}^\infty(x_0^{p-1}x_{pk+i}+p^{-1}N_{pk+i})p^{pk+i+1}\right).
\end{equation}

We have
$$
N_{pk}=\sum_{{m_0,\dots,m_{pk-1}:\atop
\sum_{i=0}^{pk-1}m_i=p,\ \
\sum_{i=1}^{pk-1}im_i=pk}} {p!\over m_0!\dots m_{pk-1}!}x_0^{m_0}\dots x_{pk-1}^{m_{pk-1}}=$$
\begin{equation}\label{e10}
=p(p-1)x_0^{p-2}x_1x_{pk-1}+\tilde{N}_{pk},
\end{equation}
where $\tilde{N}_{pk}$ does not depend on $x_{pk-1}$; moreover $p\nmid \tilde{N}_{pk}$.

Using (\ref{e10}), from (\ref{e9}) we get
$${\rm RHS \ \ of \ \ (\ref{e8})}=p^{p\gamma(x)}
\left(x_0^p+\sum_{i=1}^{p-2}(x_0^{p-1}x_i+p^{-1}N_i)p^{i+1}+\right.$$
$$\left.\sum_{k=1}^\infty (x_0^{p-1}x_{pk-1}+\tilde{N}_{pk}+p^{-1}N_{pk-1})p^{pk}
+\right.$$
$$ \sum_{k=1}^\infty (x_0^{p-1}x_{pk}-x^{p-2}_0x_1x_{pk-1})p^{pk+1}+$$ $$ \sum_{k=1}^\infty(x_0^{p-1}x_{pk+1}+x_0^{p-2}x_1x_{pk-1}+p^{-1}N_{pk+1})p^{pk+2}+$$
\begin{equation}\label{e11}
\left.
\sum_{i=2}^{p-2}\sum_{k=1}^\infty(x_0^{p-1}x_{pk+i}+p^{-1}N_{pk+i})p^{pk+i+1}\right).
\end{equation}

Consequently, $\gamma(a)=p\gamma(x)$, and $x^{p}_0\equiv a_0+a_1p\mod p^2$, $x_0=a_0$.

{\it Sufficiency.} Assume that  $a$ satisfies the conditions $(i)$ and $(ii)$. We construct a solution $x$ of the equation $x^p=a$ using the similar process of reduction to canonical form as in the proof of Theorem \ref{t2}.
First put
$$\gamma(x)={1\over p}\gamma(a).$$

Denote $x_0=a_0$ and let $M_1$ be such that $x_0^q=a_0+a_1p+M_1p^2.$

Proceeding in this way, since $x_0^{p-1}$ is not a multiple of $p$
and taking into account that integer numbers $N_k$ (see Lemma
\ref{l2}) depend only on $x_0,x_1,\dots,x_{k-1}$ and
$\tilde{N}_{pk}$ depends only on $x_0,x_1,\dots,x_{pk-2}$ we find
the existence of $x_n$ and introduce corresponding number
$M_{n+1}$ consequently for each $n\geq 1$ such that the following
congruences hold:
$$1)\,x_0^{p-1}x_i+p^{-1}N_i+M_i\equiv a_{i+1}\mod p,$$ therefore, there exists $M_{i+1}$ such that
$x_0^{p-1}x_i+p^{-1}N_i+M_i=a_{i+1}+M_{i+1}p$ for $0\leq x_i\leq p-1, \ \ i=1,\dots,p-2;$
$$2)\,x_0^{p-1}x_{pk-1}+\tilde{N}_{pk}+p^{-1}N_{pk-1}+M_{pk-1}\equiv a_{pk}\mod p,$$
therefore, there exists $M_{pk}$ such that
$x_0^{p-1}x_{pk-1}+\tilde{N}_{pk}+p^{-1}N_{pk-1}+M_{pk-1}=a_{pk}+M_{pk}p$ for $k=1,2,\dots;$
$$3)\,x_0^{p-1}x_{pk}-x^{p-2}_0x_1x_{pk-1}+M_{pk}\equiv a_{pk+1}\mod p,$$
therefore, there exists $M_{pk+1}$ such that
$x_0^{p-1}x_{pk}-x^{p-2}_0x_1x_{pk-1}+M_{pk}=a_{pk+1}+M_{pk+1}p$ for $k=1,2,\dots;$
$$4)\,x_0^{p-1}x_{pk+1}+x_0^{p-2}x_1x_{pk-1}+p^{-1}N_{pk+1}+M_{pk+1}\equiv a_{pk+2}\mod p,$$
therefore, there exists $M_{pk+2}$ such that
$x_0^{p-1}x_{pk+1}+x_0^{p-2}x_1x_{pk-1}+p^{-1}N_{pk+1}+M_{pk+1}=a_{pk+2}+M_{pk+2}p$ for $k=1,2,\dots;$
$$5)\,x_0^{p-1}x_{pk+i}+p^{-1}N_{pk+i}+M_{pk+i}\equiv a_{pk+i+1}\mod p, $$
therefore, there exists $M_{pk+i+1}$ such that
$x_0^{p-1}x_{pk+i}+p^{-1}N_{pk+i}+M_{pk+i}=a_{pk+i+1}+M_{pk+i+1}p$ for $i=2,\dots,p-2,\ \ k=1,2,\dots;$

Now making the substitutions above into equality (\ref{e11}) we
obtain
$${\rm RHS \ \ of \ \ (\ref{e11})}=p^{p\gamma(x)}
\left(a_0+a_1p+M_1p^2+\sum_{i=1}^{p-2}(a_{i+1}-M_i+M_{i+1}p)p^{i+1}+\right.$$
$$\sum_{k=1}^\infty (a_{pk}-M_{pk-1}+M_{pk}p)p^{pk}+\sum_{k=1}^\infty (a_{pk+1}-M_{pk}+M_{pk+1}p )p^{pk+1}+$$ $$\left.\sum_{k=1}^\infty(a_{pk+2}-M_{pk+1}+M_{pk+2}p)p^{pk+2}+\sum_{i=2}^{p-2}\sum_{k=1}^\infty(a_{pk+i+1}-M_{pk+i}+M_{pk+i+1}p)p^{pk+i+1}\right)=$$
$$p^{p\gamma(x)}\left(\sum_{k=1}^\infty a_kp^k +M_1p^2-\sum_{k=1}^\infty (M_k-M_{k+1}p)p^{k+1}\right)=p^{\gamma(a)}\sum_{k=1}^\infty a_kp^k=a.$$
Hence, we found the solution $\displaystyle x=\sum_{k=0}^\infty x_kp^k$ of the equation $x^p=a.$\endproof

 \begin{cor}\label{c3} Let $p$ be a prime number.

  a) The numbers $\varepsilon \in {\mathcal E}_1=\{1\}\cup\{i+jp: i^p\ \ \mbox{is not equal}\ \ i+jp \ \ \mbox{modulo} \ \ p^2\}$,  $\delta\in {\mathcal E}_2=\{p^j: j=0,\dots,p-1\}$ and products $\varepsilon\delta$ are not $p$-th power of some $p$-adic numbers.

  b) Any $p$-adic number $x$ can be represented in one of the following forms $x=\varepsilon\delta y^p$, for some $\varepsilon\in {\mathcal E}_1$, $\delta\in {\mathcal E}_2$ and $y\in \Q_p$.
  \end{cor}
\proof

a) Follows from Theorem \ref{t3}.

b)
If $x=x_0+x_1p+\dots\ne y^p$ for any $y\in \Q_p$ then by Theorem \ref{t3} we have $\varepsilon=x_0+x_1p\in {\mathcal E}_1$. We shall show that
${x\over \varepsilon}={x_0+x_1p+x_2p^2+\dots\over x_0+x_1p}=b_0+b_1p+b_2p^2+\dots$ is $p$-th power of some $y\in \Q_p$, i.e., we check the conditions of Theorem \ref{t3}: since  $\gamma(x/\varepsilon)=0$ condition $(i)$ is satisfied; we have $x_0\equiv x_0b_0 \mod p$ and $x_0+x_1p\equiv x_0b_0+(x_0b_1+x_1b_0)p\mod p^2$ which implies
 that $b_0=1$ and $b_1=0$, consequently $b_0^p\equiv b_0+b_1\mod p^2$
 i.e., the condition (ii) is satisfied.

Thus if $x\in \Q_p$ has the form $x=y^p$, then $\varepsilon=\delta=1$. If $x=x_0+x_1p+\dots$ is not $p$-th power of some $p$-adic number with $\gamma(x)=pN+j$, then we take $\varepsilon=x_0+x_1p$,
$\delta=p^j$ then $x=\varepsilon\delta y^p$ for some $y\in \Q_p$.
\endproof

Note that for a given $i\in \{1,\dots,p-1\}$ the congruence $i^p\equiv i+jp\mod p^2$ is not satisfied for some values of $j$.
For example, if $p=3$, then for $j=1$ the congruence is not true. Using computer we get the following\\

\begin{tabular}{|l|l|}
\hline
 & \mbox{Values of}\, $j$ \, \mbox{s.t.}\, $i^p\equiv i+jp\mod p^2$ \, \mbox{has no solution}\, $i\in \{1,\dots,p-1\}$\\
 \hline
p=3 & 1\\
\hline
p=5 & 2\\
\hline
p=7 & 1, 3, 5\\
\hline
p=11 & 1, 4, 5, 6, 9\\
\hline
p=13& 2, 3, 4, 8, 9, 10\\
\hline
p=17& 1, 5, 8, 11, 15\\
\hline
p=19& 4, 7, 8, 9, 10, 11, 14\\
\hline
p=23& 3, 4, 6, 9, 10, 12, 13, 16, 18, 19\\
\hline
p=29& 3, 5, 10, 11, 12, 13, 15, 17, 18, 23, 25\\
\hline
p=31& 1, 2, 5, 6, 8, 9, 11, 15, 19, 21, 22, 24, 25, 28, 29\\
\hline
p=37& 1, 4, 5, 6, 7, 10, 13, 14, 16, 20, 22, 26, 29, 30, 31, 32, 35\\
\hline
p=41& 2, 4, 6, 8, 10, 16, 24, 26, 30, 32, 34, 36, 38\\
\hline
\end{tabular}

 \begin{rk} Using this table  and Corollary \ref{c3} we obtain

 $p=3$: Any $x\in \Q_3$ has form $x=\varepsilon\delta y^3$, where $\varepsilon\in \{1,4,5\}$, $\delta\in \{1,3,9\}$;

 $p=5$: Any $x\in \Q_5$ has form $x=\varepsilon\delta y^5$, where $\varepsilon\in \{1,11,12,13,14\}$, $\delta\in \{1,5,25,125,625\}$;

 $p=7$: Any $x\in \Q_7$ has form $x=\varepsilon\delta y^7$, where $\varepsilon\in \{1,8,9,10,11,12,13,22,23,24,25,\\
 26,27,28,36,37,38,39,40,41,42\}$, $\delta\in \{7^i: i=0,1,\dots,6\}$.
\end{rk}
 3.3. {\sl THE CASE $q=mp^s$.} As is presented, the proofs of the sufficient part of the solvability criteria in Theorems \ref{t2} and \ref{t3} are given in a constructive method. Thus, not only that we know the existence of the solution, but we possess the algorithm for the construction of the solution in these cases. After cases 3.1 and 3.2 it remains the case $q=mp^s$ with some $m, s\in \N$, $(m,p)=1$. Here we shall show that this case can be reduced to cases 3.1 and 3.2: we have to find the solvability condition for $x^{mp^s}=a$. Denoting $y=x^{p^s}$, we get $y^m=a$, which is the equation considered in  case 3.1. Assume for the last equation the solvability condition is satisfied and its solution is $y={\tilde y}$. Then we have to solve
 $x^{p^s}={\tilde y}$; here we denote $z=x^{p^{s-1}}$ and get $z^p={\tilde y}$. The last equation is the equation considered in case 3.2. Suppose it has a solution $z={\tilde z}$ (i.e. the conditions of Theorem \ref{t3} are satisfied) then we get $x^{p^{s-1}}={\tilde z}$ which again can be reduced to the case 3.2. Iterating the last argument after $s-1$ times we obtain $x^p={\tilde a}$ for some ${\tilde a}$ which is also an equation corresponding to case 3.2.    Consequently,  by this argument we establish solvability condition of equation $x^{mp^s}=a$, which will be a system  of solvability of conditions for equations considered in cases 3.1 and 3.2.

\section*{ Acknowledgements}
The first author was supported by Ministerio
de Ciencia e Innovaci\'on (European FEDER support included), grant
MTM2009-14464-C02-02. The second and third authors thank the Department of Applied Mathematics, E.U.I.T. Forestal, University of Vigo,  Pontevedra, Spain,  for providing financial support
of their visit to the Department.

We thank F.M.Mukhamedov and M.Saburov for useful discussions on
previous version of our paper. Taking it to account their comments
we improved the style of proofs of Theorems 3.2 and 3.7.
{}
\end{document}